\documentclass{article}
\usepackage[all]{xy}
\usepackage[a4paper,left=3cm,top=3.5cm]{geometry}
\usepackage{verbatim}
\usepackage{CJK,CJKspace,CJKnumb}
\usepackage{latexsym,bm,easymat}
\usepackage{amsmath,amscd,amssymb,amsfonts,mathdots,ntheorem}
\usepackage{indentfirst}
\usepackage[colorlinks,linkcolor=blue,citecolor=red]{hyperref}
\usepackage{xcolor,graphicx,blkarray}
\usepackage{latexsym}
\usepackage{longtable}
\usepackage{colortab}
\usepackage{cases}
\usepackage{slashed}

\usepackage{shorttoc}
\begin{document}

%%%%%%%%%%%define of the proof=============
\def\comp{\ensuremath\mathop{\scalebox{.6}{$\circ$}}}
\def\QEDclosed{\mbox{\rule[0pt]{1.3ex}{1.3ex}}} %
\def\QEDopen{{\setlength{\fboxsep}{0pt}\setlength{\fboxrule}{0.2pt}\fbox{\rule[0pt]{0pt}{1.3ex}\rule[0pt]{1.3ex}{0pt}}}}
\def\QED{\QEDopen} %
\def\pf{\noindent{\bf Proof}} %
\def\endpf{\hspace*{\fill}~\QED\par\endtrivlist\unskip \hfill}
\def\Dirac{\slashed{D}}
%%%%%%%%%%%%%%%% Define the symbol and notation =====
\def\Hom{\mbox{Hom}}
\def\ch{\mbox{ch}}
\def\Ch{\mbox{Ch}}
\def\Tr{\mbox{Tr}}
\def\End{\mbox{End}}
\def\Re{\mbox{Re}}
\def\cs{\mbox{Cs}}
\def\spin{\mbox{spin}}
\def\Ind{\mbox{Ind}}
\def\Hom{\mbox{Hom}}
\def\ch{\mbox{ch}}
\def\cs{\mbox{cs}}
\def\nTbe{\nabla^{T,\beta,\epsilon}}
\def\nbe{\nabla^{\beta,\epsilon}}
\def\spinc{$\mbox{spin}^c$}
\def\colim{\mbox{colim}}
%%%%%%%%%%
%email-address
\newcommand{\at}{\makeatletter @\makeatother}
\renewcommand{\contentsname}{\center{Content}}
\newtheorem{defi}{Definition}[section]
\newtheorem{thm}[defi]{\textbf{Theorem}}
\newtheorem{conj}[defi]{\textbf{Conjecture}}
\newtheorem{cor}[defi]{\textbf{Corollary}}
\newtheorem{lemma}[defi]{\textbf{Lemma}}
\newtheorem{exa}[defi]{\textbf{Example}}
\newtheorem{prop}[defi]{Proposition}
\newtheorem{assump}{Assumption}
%%%%%%%%-----------------title

\title{Index of transverse Dirac operator and cohomotopy Seiberg-Witten invariant for codimension $4$ Riemannian foliation}
\date{}
\author{Dexie Lin}

\maketitle

%%%%%%%%%%%%%%%%%%%%%%%%%%%%%%%%%%%%%
\begin{abstract}
  For  closed manifolds endowed
with a Riemannian foliation of codimension $4$, one can define a transversal Seiberg-Witten map. We show that there is a finite dimensional approximation for such a map.   By such a method and  under the condition that $H^1_b(M)\cap H^1(M,\mathbb Z)$ is a lattice of $H^1_b(M)$, we can  define a foliated version of Bauer-Furuta invariant. Moreover, if the basic cohomological group is of zero dimension, we can give an estimate for the index  of transversal Dirac operator of a foliated spin structure. Furthermore, under the condition that $H^\pm_b(M)=1$, we show the vanishing of the index of the transverse Dirac operator. This gives a topological condition for the vanishing of the index of the transverse Dirac operator.
\end{abstract}
 \section{Introduction}
The interaction between geometry in dimension four is a theme which runs through a great deal of the work by many mathematicians on gauge theory over the past two decades. Moreover,
gauge theory, e.g. Seiberg--Witten equation, is one of the main tools in the study of the differential
topology of low dimensional manifolds. Since the foundational paper \cite{W} by Witten, a lot of work has
been done to apply this theory to various aspects of three and four-dimensional
manifolds. From the viewpoint of analysis, gauge theory is closely related to the study of (nonlinear)Fredholm operator and its index.  A natural idea to study Fredholm operator on foliated manifolds is to extend the framework of gauge theory to higher-dimensional situation. For instance, the high dimensional  Yang--Mills instantons
 has already been investigated by Wang \cite{Wang}, under the taut-condition the compactness of the basic Seiberg-Witten moduli space
 is showed by  Kordyukov,  Lejmi and Weber \cite{KLW}. The
theme of this article is to extending familiar constructions in gauge theory, especially Seiberg--Witten map, associated to problems in four-dimensional geometry, to higher dimensional situations, in the presence of Riemannian foliation geometric structure. Explicitly,
we consider foliation structures and develop a finite dimension approximation for smooth manifolds with foliations. Our foliation case
provides a good prototype to test and confirm a number of new characteristics about the higher dimensional Seiberg--Witten map. In particular, it is a widely generalized theory for orbifolds, which can often be realized as a particular class of foliations.

  The study
of foliations falls into largely two parts, one can study the leaf geometry or one can study the transverse elliptic operator. The leaf geometry consists of studying the individual submanifolds and how they lie within the manifold whereas the transverse geometry
is concerned with the quotient topology on the partition. On the other hand, the study of transversally elliptic operators was initiated in the seminal work of
Atiyah  \cite{Atiyah}. In summary, on a manifold with Riemannian foliation $(M,F)$,  transverse theory is to study the elliptic operator on the quotient bundle $TM/TF$.
Not long time ago,  Br\"uning et al \cite{BKR} give  a formula to express the index for transversal Dirac operator, however it is very difficult to apply the formula in reality, due to the complexity of the index formula.
Here, we will show the finite dimensional approximation for the transversal Seiberg-Witten map. One result of the paper is to define a variant Bauer-Furuta invariant with the taut condition of the foliation, using this finite dimensional approximation.
  Moreover, if  we assume that   the first basic cohomology vanishes, then there is an up bound for the index of the transversal Dirac operator(see Definition \ref{defi-transverse-Dirac} of transversal Dirac operator) of  a foliated spin structure. The  theorem below can summarize the main result of this paper.

\begin{thm}\label{thm-main-estimate}
  Let $(M,F)$ be a closed oriented manifold with a codimension $4$-foliation. If we assume that $M$ admits a foliated spin structure and the first basic de Rham cohomology vanishes, i.e. $H_b^1(M)=0$, then we have an up-bound estimate,
  \[ Ind_{\mathbb R} (\Dirac^+_b) \leq 2(b^+_b-1),\]
  where $\Dirac_b$ denotes the transversal Dirac operator and $b^+_b$ denotes the dimension of self-adjoint part of the second basic cohomology.
\end{thm}

Assuming the above theorem, one has the following estimate.

\begin{thm}
  Let $(M,F)$ be a closed oriented manifold with a codimension $4$-foliation. If we assume that $M$ admits a foliated spin structure and the first basic de Rham cohomology vanishes, i.e. $H_b^1(M)=0$, then we have an up-bound estimate,
  \[ -2(b^-_b-1) \leq Ind_{\mathbb R} (\Dirac^+_b) \leq 2(b^+_b-1),\]
  where $\Dirac_b$ denotes the transversal Dirac operator, $b^+_b$ denotes the dimension of self-adjoint part of the second basic cohomology and $b^-_b$ denotes the dimension of anti-self-adjoint part of the second basic cohomology .
\end{thm}
\begin{pf}
  It suffices to prove the lower-bound estimate.
  We reverse the orientation of this manifold by reversing the orientation of the transverse bundle $Q$. Since $Ind_{\mathbb R}(\Dirac^+_b)=-Ind_{\mathbb R}(\Dirac^-_b)$, we get that
  \[-Ind_{\mathbb R}(\Dirac^+_b)\leq 2(b^-_b-1),\]
  which proves the theorem.
\end{pf}

\noindent
For the transverse Dirac operator $\Dirac_A$, Br\"uning,  Kamber,  and Richardson gave an expression for its index \cite{BKR}. They showed that
  \[Ind(\Dirac)=\int_{\bar M_0/\bar F}A_{0,b}|\tilde{dx}|+\sum^r_{j=1}\beta(M_j),\]
  \[\beta(M_j)=\frac12\sum_\tau\frac1{n_\tau rank(W^\tau)}(-\eta(D^{S^+,\tau}_j)+h(D^{S^+,\tau}_j))\int_{\bar M_j/\bar F}A^{\tau}_{j,b}(x)|\tilde{dx}|,\]
  where the integrands $A_{0,b},~A^{\tau}_{j,b}(x)$ are similar to Atiyah-Singer integrands and notations are explained in their paper.  Due to the complexity of the formula, there is a natural question to ask:
    \vspace{2mm}

    \noindent{\bf Question}: Under what topological the Index of the transverse Dirac operator vanishes?

\vspace{2mm}
\noindent
By the theorem, we also give a lower bound for the index of the basic Dirac operator. Hence, we have  a topological condition for the vanishing of the index.
  \begin{cor}
  Let $(M,F)$ be a closed oriented manifold with a codimension $4$-foliation. If that $M$ admits a foliated spin structure and the first basic de Rham cohomology vanishes, i.e. $H_b^1(M)=0$, and $H^2_b(M)=2$, $H^+_b(M)=1$, then we have
  \[Ind_{\mathbb R}(\Dirac_b)=0.\]
\end{cor}

\noindent
The definition of the notations $H^1_b(M)$, $H^2_b(M)$ and $H^\pm_b(M)$ will be introduced in Section 2.

\noindent One  notices that our vanishing corollary does not need to establish classical Weitzenb\"ock-like formulas and the positivity on transversal scalar curvature.

 The plan of the article is the following: in Section 2, we review the basic notions and properties of the Riemannian foliation; in Section 3, we  review some properties of Fredholm, stable cohomotopy and cohomotopy group and give an explicit construction of the high dimensional Bauer-Furuta invariant.; in Section 4, we give a proof of Theorem \ref{thm-main-estimate}.

\vspace{0.5cm}

{\bf Acknowledgement:} The author thanks Mikio Furuta for the discussion and encouragement to study gauge theory.

 \section{Geometry of Foliation }

In this section, we review some results of the previous work. In the first subsection, we review the classical results of the geometric foliation.
Let $M$ be a closed oriented $n$ dimensional manifold with dimension $p$ oriented foliation $F$. We  denote codimension of this foliation by $q=n-p$. For more details of this  section, we give a reference \cite{Ton2}.
 \begin{defi}
   A Riemannian metric $g_Q$ on $Q$ is said to be bundle-like, if
   \[L_Xg_Q\equiv0,\]
   for any $X\in \Gamma(F)$, where $Q=TM/TF$. We say $(M,F)$ is a Riemannian foliation, if $Q=TM/TF$ admits a bundle-like metric.
 \end{defi}

 Given a metric $g$ on $TM$, $Q$ can be identified with the orthogonal complement to $TF^\perp$ by $g$ and inherits  a metric $g_{F^\perp}$. We have the following equivalence,
 \[\mbox{a metric }g\mbox{ of }TM\mbox{ corresponds a triple }(g_F,\pi_F, g_Q),\]
 where $g_F=g|_{TF}$ and $\pi_F$ is the projection $TM\to TF$.

\noindent
   Let $M$ be a manifold with foliation $F$.  A Riemannian metric $g$ on $TM$ is said \emph{bundle-like}, if the induced metric $g_Q$ is bundle-like.
 By the work of Reinhart $\cite{Reinhart}$, it is known  that the bundle-like metric can be locally written as $g_Q=%\sum_{i,j}g_{ij}(x,y)\omega^i\otimes \omega^j+
 \sum_{k,l}g_{k,l}(y)dy^k\otimes dy^l$, where $(x,y)$ is in the foliated chart of $M$.% and $\omega^i=dx^i+a^i_\alpha(x,y) dy^\alpha$.
In this paper, we always assume that $(M,F)$ is a Riemannian foliation.

Let $\pi$ be the canonical projection $TM\to Q$.  We define a connection $\nabla^{T}$ on $Q$, by
$$\nabla^{T}_Xs=\begin{cases}
   \pi([X,Z_s])& X\in \Gamma(F),\\
   \pi(\nabla_X Z_s)& X\in \Gamma(F^\perp),
 \end{cases}$$
 for any section $s\in\Gamma(Q)$,
 where $Z_s\in \Gamma(TM)$ is a lift of $s$, i.e. $\pi(Z_s)=s$. If $(M,F)$ is a Riemannian foliation, then by the Koszul-formula \cite[Theorem 5.9]{Ton2},  we have that $\nabla^T$ is  uniquely determined by $g_Q$. Moreover, one can verify that it is  torsion free and metric-compatible, whose leafwise restriction coincides with the Bott-connection. We set $R^T$ as the curvature of this connection.
 We define the transverse Ricci curvature and scalar curvature by
 \[Ric^{T}(Y)=\sum^{q}_{i=1}R^{T}(Y,e_i)e_i,~
 Scal^{T}=\sum^{q}_{i=1}g_Q(Ric^{T}(e_i),e_i),\]
 where $\{e_i\}$ is a local  orthonormal frame of $Q$.
 % One has that $R^T$ satisfies the   condition \[\iota_XR^T=0,\]for any $X\in\Gamma(F)$.
  We define the basic forms as follows:
 \[\Omega^r_b(M)=\{\omega\in\Omega^r(M)\big|~\iota_X(\omega)=0,~L_X(\omega)=0,\forall X\in\Gamma(F)\}.\]

 By the work of Alvarez L\`opez \cite{AL}, we have the following $L^2$ orthogonal decomposition  for the forms on $M$,
 \[\Omega(M)=\Omega_b(M)\oplus\Omega_0(M).\]
% The $L^2$ -norm is defined below. We  have the basic Hodge-star operator, \[\bar*:\bigwedge^rQ^*\to \bigwedge^{q-r}Q^*.\]
% Choosing  a local orthonormal basis $\{e_i\}_{1\leq i\leq p}$ of $TF$, we define the character form of the foliation $\chi_F$ by,
% \[\chi_F(Y_1,\cdots, Y_p)=\det(g_F(e_i,Y_j))_{1\leq i,j\leq p},\]
% for any section $Y_1,\cdots, Y_p\in \Gamma(TM)$.
 The basic Hodge-star operator is  related to the usual Hodge-star operator by the formula $\bar*\alpha=(-1)^{(q-r)\dim(F)}*(\alpha\wedge\chi_F)$, where $\chi_F$  is the characteristic form of the foliation \cite[Formula 4.16]{Ton2}
 We have  the volume density formula,
 $dvol_M=dvol_Q\wedge\chi_F$.
 For the section $\alpha\in \Gamma(\bigwedge^rQ^*)$, we define its $L^2$ norm by
 $\|\alpha\|^2_{L^2}=\int_M\alpha\wedge\bar*\alpha\wedge\chi_F$.
% For any section $\alpha\in \Gamma(\bigwedge^rQ^*)$, we have that \[\int_M\alpha\wedge\bar*\alpha\wedge\chi_F=\int_M\alpha\wedge*\alpha.\]
For the bundle-like metric, we  have $\bar*:\Omega^r_b(M)\to\Omega^{q-r}_b(M).$

\begin{defi}
  The mean curvature vector field is defined by $\tau=\sum^{\dim F}_{i=1}\pi(\nabla_{\xi_i}\xi_i)\Gamma(Q)$, where $\{\xi_i\}$ is a local orthonormal basis of $F$. Let $\kappa\in \Gamma(Q^*)$ be the dual to $\tau$ via the metric  $g_Q$.
\end{defi}

%\begin{prop}[Rummler \cite{Rummler}]\label{formula-Rummler}For any metric $g$ on $TM$,  we get  \[d\chi_F=-\kappa\wedge\chi_F+\phi_0,\]  where $\phi_0$ belongs to $F^2\Omega^{p+1}=\{\omega\in\Omega^{p+1}(M)\big|\iota_{X_1}\cdots\iota_{X_p}\omega=0,\mbox{ for any }X_1,\cdots,X_p\in \Gamma(F)\}$. This, implies that $w\wedge\phi_0=0$ for any $w\in \Gamma(\bigwedge^{q-1}Q^*)$.\end{prop}

% By direct calculation \cite[Chpater 7]{Ton2},    we have that the $L^2$ formal adjoint  of $d$ is $\delta=(-1)^{m(*+1)+1}\bar*(d-\kappa\wedge)\bar*$.
% \end{prop}
% \begin{pf}   By direct calculation, we have for any sections $\alpha,~\beta\in \Omega^*(Q^*)$, the following identity holds:   \begin{eqnarray*}     d(\alpha\wedge\bar*\beta\wedge\chi_F)&=&((d\alpha)\wedge\bar*\beta)\wedge\chi_F+   (-1)^{r-1}(\alpha\wedge(d\bar*\beta))\wedge\chi_F\\   &&+(-1)^{m-1}(\alpha\wedge\bar*\beta)\wedge   (-\kappa\wedge\chi_F)+(\alpha\wedge\bar*\beta)\wedge\phi_0.   \end{eqnarray*}   When $\alpha\in\Omega^{r-1}(Q^*)$ and $\beta\in\Omega^r(Q^*)$, one can deduce that $(\alpha\wedge\bar*\beta)\wedge\phi_0\equiv0$. Taking the integral of the above formula, we have   \begin{eqnarray*}     (d\alpha,\beta)_{L^2}&=&(-1)^r\int_M\alpha\wedge(d\bar*\beta-(-1)^{m-r}    \bar*\beta\wedge\kappa)\wedge\chi_F\\     &=&(-1)^{m(r+1)+1}\int_M\alpha\wedge\bar*(\bar*(d-\kappa\wedge)\bar*\beta)\wedge     \chi_F.   \end{eqnarray*} \end{pf}

By the decomposition, we have that  $\kappa=\kappa_b+\kappa_0$, where $(\kappa_0,\omega_b)_{L^2}=0$ for any basic one form $\omega_b$.
We call $\kappa_b$ the basic mean curvature form.  It is known that
$d\kappa_b=0$
 and the cohomology class $[\kappa_b]$ is independent on any bundle-like metric \cite{AL}.

 \begin{defi}
  We say a foliation is \emph{taut}, if  there is a metric on $M$ such that $\kappa=0$, i.e. all leaves are minimal submanifolds.
\end{defi}

For a fixed Riemannian foliation $F$, the taut condition has a topological obstruction.

\begin{prop}[ Tondeur {\cite[Page 96]{Ton2}}]\label{prop-taut}
  Let $(M,F)$ be a Riemannian foliation. Suppose that $M$ is closed oriented and each leaf is also oriented. Then, the following statements are equivalent
  \begin{itemize}
    \item $H^q_b(M)\neq0$, $q$ is the codimension of this foliation $F$.
    \item $[\kappa_b]=0~\in H^1_b(M)$,
    \item the foliation  admits a taut metric.
  \end{itemize}
\end{prop}

%\begin{prop}[ Dominguez \cite{D}]  Any Riemannian foliation $F$ carriestense bundle-like metrics, i.e. having basic mean curvature form\[\kappa=\kappa_b.\]\end{prop}\textbf{Remark}: It is known that any bundle-like metric can be deformed inthe leaf directions leaving the transverse part unchanged in such a way that the mean curvature form becomes basic.In this article, we always let  $\kappa$  be basic, i.e. $\kappa=\kappa_b$.

 \begin{prop}[Tondeur {\cite[Theorem 7.18]{Ton2}}]
   Let $d_b$ denote the restriction of $d$ on the basic forms. %For any basic form $\alpha$ we define the $L^2$-norm by \[\|\alpha\|^2_{L^2}=\int_M\alpha\wedge\bar*\alpha\wedge\chi_F.\]
   Then the  $L^2$-formal adjoint   of $d_b$ is $\delta_b=(-1)^{m(*+1)+1}\bar*(d_b-\kappa_b\wedge)\bar*$.
 \end{prop}
 We define the basic Laplacian operator by $\Delta_b=d_b\delta_b+\delta_bd_b$.

In order to define the basic Dirac operator, we need to review the definition of foliated bundle and the associated notions.

 \begin{defi}
   A principal bundle $P\to M$ is called foliated, if it is quipped with a lifted foliation $F_P$ invariant under the structure group action, transversal to  the tangent space to the fiber and $TF_P$ projects isomorphically onto $TF$. We say a vector bundle $E\to M$ is foliated, if its principal bundle $P_E$ is foliated.
 \end{defi}

 \begin{defi}
   A connection $\omega$ of the foliated principal bundle $P$ is called adapted, if the horizontal distribution associated to this connection   contains the foliation $F_P$. A covariant connection on a foliated vector bundle is called adapted, if its associated connection on the principal bundle is. An adapted connection $\omega$ is called basic, it is a Lie algebra valued basic form. The similar notion for the covariant connection.
 \end{defi}

% For a foliated vector bundle $E$ over $M$, by the definition, it is  known that for any two connections $\nabla^1$ and $\nabla^2$ adapted to this foliated vector bundle, we have that \[\nabla^1_Vs=\nabla^2_Vs,\] for all $s\in\Gamma(M,E)$ and $V\in\Gamma(F)$.
 Using a adapted connection,  we define the basic sections by \[\Gamma_b(E)=\{s\in \Gamma(E)\big|~ \nabla_Xs\equiv0,~\mbox{for all }X\in \Gamma(F)\},\]
 where $\nabla$ ia an adapted connection. It is  known that the space of  basic sections is independent of the choice of the adapted connection.
 \begin{defi}
   A transverse Clifford module $E$ is a complex vector bundle over $M$ equipped with a hermitian metric  satisfying the following properties.
 \begin{enumerate}
   \item  $E$ is a  bundle of $Cl(Q)$-modules, and the Clifford action $Cl(Q)$ on $E$ is skew-symmetry, i.e.  \[(s\cdot\psi_1,\psi_2)+(\psi_1,s\cdot\psi_2)=0,\]
 for any $s\in\Gamma(Q)$ and $\psi_1,\psi_2\in\Gamma(E)$;
   \item $E$ admits a basic metric-compatible connection, and this connection is compatible with the Clifford action.
 \end{enumerate}
 \end{defi}

 We say $(M,F)$ admits a transverse \spinc structure, if $Q$ is a foliated \spinc structure. And the \spinc structure corresponds to the foliated line bundle admitting a basic connection.

% {\bf Remark:} Actually, for any  complex Hermitian foliated bundle $E$,   such that $E$ is a transverse Clifford module, self-adjoint and equipped with a Hermitian Clifford connection $\nabla^E$, by \cite{KT} we can choose $\nabla^E$ to be a \emph{basic connection} which means that theconnection and curvature forms of $\nabla^E$ are  basic forms. In this paper, we always assume that the connection is basic.

%\hfill

% For a foliated \spinc principle bundle $P$ over $M$, the automorphism  group $G_b$ is a subgroup of $G=\{u:M\to U(1)\}$  satisfies that\[Vu=0,\] for all $V\in\Gamma(F)$ for $u\in G$.

% {\bf Example}:For $T^4$ with dense a $1$-dimensional foliation $V$.  We have that for any foliated \spinc structure, the automorphism group $G_b$ consists of all constant maps.To see this, we fix a point $x_0\in T^4$, choose a metric on $T^4$. For any $x\in T^4$, and $\epsilon>0$, there is a geodesic ball at $x_0$ with radius $\delta$, such that for any point $x$ in this ball, we have \[|f(x)-f(x_0)|<\epsilon,\]for any smooth function $f$. As the foliation is dense, there is a constant $T(\delta)$, such that $\exp(T(\delta)x_0)$ locates in this ball. Since for any $u\in G_b$, $u$ is constant along this orbit $\exp(tx_0)$, by the above argument, we get \[|u(x)-u(x_0)|<\epsilon,\] hence $u(x)=u(x_0)$.

\begin{defi}\label{defi-transverse-Dirac}
   Fixing a basic connection $\nabla^E$, we define the Dirac operator $\Dirac_b$ by $\Dirac_b=\sum^q_{i=1}e_i\cdot\nabla^E_{e_i}$ action on $\Gamma(E)$, where $\{e_i\}$ is a local orthonormal basis of $Q$.
\end{defi}
 Note that it is not formally self adjoint in general, whose adjoint operator is $\Dirac^*_b=\Dirac_b-\tau_b$.
 %For the non-taut case, we define $\Dirac_b=\sum^q_{i=1}e_i\cdot\nabla^E_{e_i}-\frac12\tau_b$.
%By straightforward calculation, we have that the transversal Dirac operator maps the basic sections $\Gamma_b(S)=\{s\in\Gamma(M,S)\big|\nabla_Xs\equiv0,\mbox{ for any }X\in\Gamma(F)\}$ to basic sections.
%   \begin{prop}   The transversal Dirac operator maps the basic sections $\Gamma_b(S)=\{s\in\Gamma(M,S)\big|\nabla_Xs\equiv0,\mbox{ for any }X\in\Gamma(F)\}$ to basic sections. \end{prop}
% \begin{pf}   By the straightforward calculation, one gets that  \begin{eqnarray*}     \nabla_X(\Dirac_b\psi)&=&\sum^m_{i=1}\nabla_X(e^i\nabla^E_{e_i}\psi)\\     &=&\sum^m_{i=1}(\nabla^T_X(e^i))\nabla^E_{e_i}\psi+     \sum^m_{i=1}e^i(\nabla^E_X\nabla^E_{e_i}\psi)\\     &=&\sum^m_{i=1}(\nabla^T_X(e^i))\nabla^E_{e_i}\psi+     \sum^m_{i=1}e^i(\nabla^E_{e_i}\nabla^E_X\psi)\\     &&+\sum^m_{i=1}e^i(\nabla_X(\nabla^E))     (e_i,\psi)+\sum^m_{i=1}e^i\nabla^E_{[X,e^i]}\psi.   \end{eqnarray*}   Letting $X\in\Gamma(F)$ and $\psi\in \Gamma_b(S)$, on the above formula we have that the second and third terms vanish,  the first and last terms cancel. \end{pf}

Let $E$ be a foliated vector bundle on $M$ equipped with a basic Hermitian structure
and a compatible basic connection $\nabla^E$, we define
\[\|u\|_{p,k}=\sum^k_{j=1}(\int_M|(\nabla^E)u^j|^pdvol_M)^{\frac1p},\]
for any $u\in \Gamma_b(E)$. Let $L^p_k$ be the complete space of $\Gamma_b(E)$ with respect to such a norm.
One has the similar Sobolev embedding and Sobolev multiplication properties for basic sections, which are shown in \cite[Theorem 9, 10, 11]{KLW}. To make this paper complete, here we just give the statements.

\begin{thm}Suppose that $(M,F)$ is a closed oriented manifold with codimension-$m$ foliation $F$, then the following inclusions hold:\begin{itemize}  \item \[L^p_k\hookrightarrow L^q_l,\]  where integers $l,~k$ satisfying $0\leq l\leq k$ and $l-\frac mk\leq k-\frac mp$.  \item \[L^p_k\hookrightarrow C^l,\]  where $l<k-\frac mp$.\end{itemize}\end{thm}

\begin{thm}%[ \cite{KA}]
Let $0\leq l\leq k$, under the setting of above theorem, we have the following continuous maps: \begin{itemize}    \item \[L^p_k\times L^q_l\to L^q_l,\]    where $k - \frac mp > 0$ and $k -\frac mp > l- \frac mq$.    In particular, if $k=l,~p=q$ and $k-\frac mp>0$, then    \[L^p_k\times L^p_k\to L^p_k.\]    \item \[L^p_k\times L^q_l\to L^r_t,\]    where $k-\frac mp<0$, $l-\frac mq<0$ and $l$ satisfies $0\leq t\leq l$,$r$ satisfies $0<\frac tm+\frac1p-\frac km+\frac1q-\frac lm\leq \frac1r\leq1$.    \item \[(L^p_k\cap L^\infty)\times(L^q_l\cap L^\infty)\to (L^q_l\cap L^\infty),\]        where $k=\frac mp$ and $l-\frac mq\leq0$.    \item\[(L^p_k\cap L^\infty)\times L^q_l \to L^q_l,\]    where $l-\frac mp<0$. \end{itemize}\end{thm}

%One can follow the idea of \cite{KA} to give a proof. Here we give another version of the proof.
%\begin{pf}  Let $(E,\nabla^E, h_E)$ be a foliated vector bundle over $(M,F)$ with basic connection  $\nabla^E$ and hermitian metric $h_E$. By \cite[Chapter 4]{Mol} or \cite{BKR}, one can find a quintuple $(W,G,E',\nabla',h_G)$, where $G$ is compact Lie group, $W$ is a closed $G$-manifold, $E'$ is a $G$-equivariant vector bundle   with a $G$-equivariant connection $\nabla'$ and $h_G$ is a hermitian metric of $E'$, such that the following properties are hold:  \begin{itemize}    \item   $\Omega^*_b(M)\cong \Omega^{*,G}(W)$, where $\Omega^{*,G}(W)$ denotes the space of $G$-invariant forms of $W$;    \item there is an isometric identification $\Gamma_b(M,E)\cong \Gamma^G(W,E')$ with respect to $h_E$ and $h_G$, where $\Gamma^G(W,E')$ denotes the space of $G$-invariant sections.    \item $Q\cong  T_GW$, where $T_GW\subset TW$ denotes the sub-bundle of $TW$, which are orthogonal to the $G$-orbit.  \end{itemize}  Since the Sobolev embeddings hold for the $G$-invariant section of $W$, by the above identification, the above two theorems can be proved.\end{pf}

\section{Codimension $4$ Seiberg-Witten and high dimensional Bauer-Furuta invariant}

In this section, we will review the work of Bauer and Furuta \cite{BF} and generalize their results to the manifolds with  foliation.

\subsection{Review of  stable cohomotopy and cohomotopy group}
In this subsection, we review some results for a general Fredholm operator with a certain condition.
 Let $H'$ and $H$ be  separable Hilbert spaces and
  \[f:H'\to H\]
be a Fredholm map, which is a sum of compact perturbation and linear Fredholm operator, i.e. $f=l+c$, where $l$ is linear Fredholm and the continuous map $c$ maps bounded sets to
subsets of compact sets.

\begin{lemma}[Bauer and Furuta {\cite[Lemma 2.2]{BF}}]
  Let $l:H'\to H$ be a continuous linear Fredholm map between Hilbert spaces, and $c:H'\to H$ be a compact map. Then, the restriction of the map $f=l+c$ to any closed and bounded subset $A'\subset H'$ is proper, i.e. $f\mid_{A'}$ is proper.
  If the preimages of bounded sets in $H$ are bounded, then $f$ is proper and extends to a proper map $$f^+:(H')^+\to H^+$$
  between the one-point completed Hilbert spheres, where $H^+$($(H')^+$) denotes the one-point completed Hilbert spheres of $H$($H'$).
\end{lemma}

 %\begin{pf}   Let $\rho :H'\to \ker(l)$ be the orthogonal projection, hence we can $f\mid_{A'}$ factors through three parts:   \[\begin{array}{rcl}     A'&\to& H\times \overline{c(A')}\times \overline{\rho(A')}\\     a'&\mapsto &(l(a'),c(a'),\rho(a'));   \end{array}\]    \[\begin{array}{rcl}     H\times \overline{c(A')}\times \overline{\rho(A')}&\to& H\times \overline{c(A')}\times \overline{\rho(A')}\\(h,s,e)&\mapsto &(h+s,s,e);\end{array}\]    \[\begin{array}{rcl}     H\times \overline{c(A')}\times \overline{\rho(A')}&\to& H\\    (h,s,e)&\mapsto &h.   \end{array}\]   Since the external factors are compact, we have $f\mid_{A'}$ is proper.   For the bounded condition, the preimages of points in $H$ are bounded and closed, hence the above argument can be also applied, which implies that the preimage of point is compact.   We show the properness of $f$. Let $h\in H$ such that $h\in \overline{f(A')}$ and $A'\hookrightarrow H'$ is a closed subset. By the bounded condition, we get $h\in \overline{f(A'_0)}$ for some $A'_0$ which is a bounded closed subset of $A'$. We know that $f\mid_{A'_0}$ is a proper map, we have $h\in f(A'_0)\subset f(A')$.   Finally, if $A'\subset (H')^+$ which is  closed, we know that $\overline{f^+(A')}$ contains the infinity point, then $f^+(A')\cap H$ is unbounded, so $A'$ is unbounded.  \end{pf}

% {\bf Remark}: Properness $\not\Rightarrow$ the boundedness, e.g. $f(x)=x+\sum_{n=1}(n-1)\phi(x-ne_n)e_n$, here $\phi:H\to [0,1]$ is supported in $[0,\frac12]$.

 We assume that $f$ satisfies the bounded condition.
 Let $W\subset H$ be a finite dimensional linear subspace and let $W'=l^{-1}(W)$ be the preimage under the   map $l$. We have that  the inclusion  $W^+\to H^+\setminus S(W^\perp)$ is a deformation retract, where $S(W^\perp)$  denotes the unit sphere in the orthogonal complement $W^\perp$ of $W$.  The retract map $\rho_W$ can described as follows:
 \begin{itemize}
   \item We identify $H^+\cong S(\mathbb R\oplus W\oplus W^\perp)$, by $h\mapsto (|h|^2+1)^{-1}(|h|^2-1,2h)$.
   \item $W^+$ maps to the equatorial $S(\mathbb R\oplus W\oplus 0)$ and $S(W^\perp)$ maps to $S(0\oplus0\oplus W^\perp)$.
   \item The retract homotopy shrinks the latitudes in $S(\mathbb R\oplus W\oplus W^\perp)\setminus S(0\oplus0\oplus W^\perp)$ to the equator.
 \end{itemize}

 {\bf Remark:}

 $\rho_W$ has the following  property : For $h\in H\setminus W^\perp$, the vector $\rho_W$ differs the projection $pr_W(h)$ to $W$ by a positive scalar, i.e.
  $\rho_W(h)=\lambda(h)pr_W(h)$.

  \begin{lemma}[Bauer and Furuta {\cite[Lemma 2.3]{BF}}]
    There exists $V\subset H$ a linear subspace, such that the following statements hold:
    \begin{itemize}
      \item[(1)] $V+Im(l)=H$.
      \item[(2)] for any $W\supset V$ with $W=U\perp V$ such that $f\mid_{(W')^+}:(W')^+=(l^{-1}(W))^+\to H^+$, such that $Im(f\mid_{(W')^+})\cap S(W^\perp)=\emptyset$.
      \item[(3)] $\rho_Wf\mid_{(W')^+}$ and $Id_{U^+}\wedge \rho_Vf\mid_{(V')^+}$ are homotopy as pointed map
          \[(W')^+\cong U^+\wedge (V')^+\to U^+\wedge V^+=W^+.\]
    \end{itemize}
  \end{lemma}

  We denote by $\pi^n(X)=[X,S^n]$ the homotopy classes of continuous maps for a CW-complex  $X$. The following definition was given by Bauer and Furuta.
%  \begin{prop}    If $X$ is of dimension $2n-2$, then $\pi^n(X)$ is an Abelian group.  \end{prop}  \begin{pf}    For any $[\alpha],[\beta]\in \pi^n(X)$, we define $[\alpha]+[\beta]$ as follows:    Let $\Delta:X\to X\times X$ as diagonal, then $(\alpha\times\beta)\comp\Delta:X\to S^n\times S^n$, for the dimension reason, there exists $f:X\to S^n\vee S^n$, such that the diagram commutes.    \[\xymatrix{  X \ar[r]^{f} \ar@{->}[dr]_{(\alpha\times\beta)\comp\Delta}                &S^n\vee S^n\ar@{^{(}-_>}[d] \\                & S^n\times S^n    }\]   Denote by $\theta:S^n\vee S^n\to S^n$ the folding map, we have that  $[\alpha]+[\beta]:=[\theta\comp f]$ is an Abelian sum.  \end{pf}
% If $(X,A)$ is a pair of CW-complex, with $\dim\leq2n-2$, then we have \[..\to\pi^i(X)\to \pi^i(A)\to \pi^i(X,A)\to \pi^{i+1}(X)\to..,\] where $A$ is a closed set and $\pi^i(X,A):=\pi^i(X/A)$. We have $\pi^n(X)\cong \pi^{n+1}(S^1\wedge X)$, with $\dim(X)\leq 2n-2$.

% We denote the stable homotopy classes  of maps by   \[\pi^n_{st}(X)=\lim_{N\to\infty}[S^N\wedge X,S^{n+N}].\]
 \begin{defi}[Bauer and Furuta~\cite{BF}]
   For a nonlinear Fredholm map $f$ satisfying the above assumption, we define a class $[f]$ by,
   $$[f]=\colim_{V\subset H}[(f\mid_{l^{-1}(V)})^+]\in \colim_{V\subset H}
   [(l^{-1}(V))^+,H^+\setminus S(V^\perp)].$$

 \end{defi}
Combine the homotopy equivalences $V^+\subset (H\setminus (S(V^\perp)))$ with an isomorphism
   \[\pi^{st}_{Ind(l)}(S^0)=\colim_{V\subset H}
   [(l^{-1}(V))^+,V^+]
   \cong \colim_{V\subset H}
   [(l^{-1}(V))^+,H^+\setminus S(V^\perp)],\] we regard $[f]$ as an element of the stable homotopy group $\pi^{st}_{Ind(l)}(S^0)$.
 Let $f:E'\to E$ be a continuous map between Hilbert bundles $E'$ and $E$ over a compact manifold
 $Y$, such that $E=Y\times H$. We extend the map $pr_H\comp f:E'\to H$ to a point compactification,
 \[(pr_H\comp f)^+:T(E')\to H\]
 where $T(E')$  denotes for the Thom space. Let $\lambda=F_0-F_1$ be  an element of the $K$-group over $Y$, such that $F_1=Y\times V$, we set
 \[\pi^n_H(Y;\lambda)=\colim_{U\subset V^\perp}[U^+\wedge TF_0,U^+\wedge V^+\wedge S^n]
 =\colim_{W\subset H}[W^+T\lambda,W^+\wedge S^n].\]
 For a compact Lie group $G$, let $H$ be a Hilbert space equip with an orthogonal $G$-action such that $H$ contains the trivial representation. We assume that the space of equivariant morphisms $Hom_G(V,H)$ for a real $G$-module $V$ is either zero or of infinite dimension. Let $\lambda=F_0-F_1$ be a virtual equivariant vector bundle over a finite $G$-CW complex $Y$ such that $F_1\cong Y\times V$ is a trivial bundle with $V\subset H$ a finite dimensional $G$-subrepresentation.
 We define the stable equivariant cohomotopy groups by
 \[\pi^n_{G,H}(Y;\lambda)=\colim_{U\subset V^\perp}[U^+\wedge TF_0,U^+\wedge V^+\wedge S^n]^G
 =\colim_{W\subset H}[W^+T\lambda,W^+\wedge S^n]^G.\]
 Summarizing the above arguments, we have the following theorem.
 \begin{thm}[Bauer-Furuta {\cite[Theorem 2.6]{BF}}]\label{thm-2.6}
   An equivariant Fredholm map $f=l+c:E'\to E$ between $G$-Hilbert bundle over $Y$ with $E\cong Y\times H$, which extends continuously to the fiberwise one-point completions, defines a stable cohomotopy class
   \[[f]\in \pi^0_{G,H}(Y;ind(l)),\]
   which is independent of the presentation of $f$ as a sum.
 \end{thm}

\subsection{High dimensional Bauer-Furuta invariant}

 In this section,   $M$ is a closed oriented taut manifold with codimension $4$ Riemannian foliation satisfying the following assumption.
 \begin{assump}\label{assum-main-1}
  Let $(M,F)$ be a oriented closed manifold with codimension 4 oriented Riemannian foliation $F$ and admits a transverse \spinc structure $\mathfrak s$. Suppose  that  $H^1_b(M)\cap H^1(M,\mathbb Z)\subset H^1(M)$ is a lattice of $H^1_b(M)$ and $(M,F)$ admits a taut bundle-like metric.
\end{assump}

   We choose a bundle-like metric $g:=g_F\oplus g_Q$.
Suppose that $Q$ admits a foliated \spinc structure, then it admits a spinor bundle $S=S^+\oplus S^-$. We get a decomposition
 \[\Lambda^2Q^*=\Lambda^+Q^*\oplus\Lambda^-Q^*,\]
 where $\Lambda^\pm$ corresponds to the $\pm1$-eigenvalue under the action of the basic Hodge operator $\bar*$. Let $\mathcal A_b$ denotes the space of basic connections of the determinant line bundle of the spinor bundle. We define the {\emph basic Seiberg- Witten} equations, by
 \begin{equation}
   \begin{cases}
   \Dirac^A_b\psi=0,\\
   F^+_A=q(\psi),
 \end{cases}\label{eqn-SW}
 \end{equation}
 for the pair $(A,\psi)\in\mathcal A_b ,~\psi\in\Gamma_b(S^+)$, where $q(\Phi)=\Phi\otimes\Phi^*-\frac{|\Phi|^2}21$ and we used the identification
\[cl_+:\Lambda^{2,+}_b\otimes\mathbb C\to End^0_b(S^+),\]
between the self-adjoint basic two forms and traceless basic  endomorphism  of $S^+$(see \cite[Chaper 3,4]{Morgan} for more details).
The moduli space $\mathcal M(\mathfrak s)$ is the space consisting  of all solutions $(A,\Phi)$ mod  the   gauge group $\mathcal G_b=C^\infty_b(M,S^1)$, where $C^\infty_b(M,S^1)$ denotes the set of $S^1$-valued basic functions of $M$ and  the gauge action is defined as following: for each $g\in\mathcal G_b$,
\[g(A,\Phi)=(A-2g^{-1}dg,g\cdot \Phi),\]
 where $\cdot$ denotes the Clifford multiplication.
% The basic gauge group $\mathcal G_b$ consists of leaf-invariant $U(1)$-valued function.
  We define the operators $l,~c:i\Omega^1_b\times\Gamma_b(S^+)\to i\Omega^+_b\times\Gamma_b(S^-)\times i\Omega^0$, by
\[\begin{array}{l}
  l(a,\Phi)=(d^+_ba,\Dirac_b\Phi,\delta_ba),\\
  c(a,\Phi)=(-\sigma(\Phi),\frac a2\cdot\Phi,0).\\
\end{array}\]

 \begin{lemma}
  The map $l:L^2_{k+1}(i\Omega^1_b(M)\oplus \Gamma_b(S^+))\to L^2_k(i\Omega^+_b(M)\oplus \Gamma_b(S^-)\oplus i\Omega^0_b(M))$ is a Fredholm operator.
\end{lemma}

\begin{pf}
  It suffices to show the Fredholmness for  operators $\Dirac_b$ and $d^+_b+\delta_b$ independently. The routine is similar to show that the elliptic operators on compact space are Fredholm.
 Let $\Dirac_b : L^2_{k+1}(\Gamma_b(S^+))\to \L^2_k(\Gamma_b(S^-))$ be the basic elliptic operator.
   By the estimate \cite[Theorem 12]{KLW}, we have that
   \[\|\psi\|_{L^2_{k+1}}\leq C(\|\Dirac_b\psi\|_{L^2_k}+ \|\psi\|_{L^2_k}),\]
   for any $\psi \in L^2_{k+1}(\Gamma_b(S^+))$, where $C$ is some positive constant depending on $k$. 
   We claim that the kernel $\ker(\Dirac_b)$ is of finite dimension. Otherwise, we can find a sequence of solutions $\{\psi_i\}_{i\geq1}$, such that they are mutually orthogonal to each other under the $L^2$-metric and   $\|\psi_i\|_{L^2}=1$ for each $i$. However, by the compact Rellich embedding, we have that $\{\psi_i\}_{i\geq1}$ must have a convergent subsequence in $L^2_k$, which is a contradiction.
   We show that the range of $\Dirac_b$ is closed. Let $\{\psi_i\}_{\i\geq1}$ be a sequence in $L^2_{k+1}(\Gamma_b(S^+))$ with $\phi_i=\Dirac_b\psi_i\to \phi\in \L^2_k(\Gamma_b(S^-))$. Without loss of generality, we assume that $\psi_i\not\in\ker(\Dirac_b)$ for each $i$.  It suffices to show that $\{\psi_i\}$ has a convergent subsequence. If $\{\psi_n\}$ is bounded, then by  the compact inclusion $L^2_{k+1}\to L^2_k$ and the estimate $$\|\psi_i-\psi_j\|_{L^2_{k+1}}\leq C(\|\phi_i-\phi_j\|_{L^2_k}+\|\psi_i-\psi_j\|_{L^2_k}),$$    we have that $\{\psi_i\}$ has a convergent subsequence.
    Assume that $\{\psi_n\}$ is not bounded. We set $\psi'_i=\psi_i/\|\psi_i\|_{L^2_{k+1}}$. It is clear that $\Dirac_b\psi'_i\to 0$. By applying the  arguments of the bounded case, we have that $$\psi'_i\to\psi',~\|\psi'\|_{L^2_{k+1}}=1  \mbox{ and } \Dirac_b\psi'=0,$$ which contradicts to the choice $\psi_i\not\in\ker(\Dirac_b)$. It remains to show that the cokernel of $\Dirac_b$ is of finite dimension, which follows from the formal self-adjoint property of $\Dirac_b$.
   Similarly, we can show the Fredholmness of $d^+_b\oplus\delta_b$.
\end{pf}

 \begin{prop}[Glazebrook and   Kamber \cite{GK}]
   Let $(M,F,g)$ be a manifold with Riemannian foliation  admitting  a transverse \spinc structure $\mathfrak s$, for the associated basic Dirac operator we have
   \[(\Dirac^A_b)^*\Dirac^A_b\psi=(\nabla^A_b)^*\nabla^A_b\psi
   +\frac14( Scal^T-\delta_b\kappa_b+|\kappa_b|^2)\psi
   +\frac12F^+_A\cdot\psi,\]
   where $(\nabla^A_b)^*\nabla^A_b=-\sum^m_{i=1}\nabla^A_{e_i}\nabla^A_{e_i}+\sum^m_{i=1}
   \nabla^A_{\nabla^T_{e_i}e_i}+\nabla^A_{\tau_b}$.
 \end{prop}

% \begin{pf}   We have   \[(\Dirac^A_b)^*\Dirac^A_b\psi=\sum^m_{i,j=1}e_i\nabla^T_{e_i}e_j\nabla^T_{e_j}\psi   -\tau_b\Dirac^A_b\psi,\]   and for the term   \begin{eqnarray*}     \sum^m_{i,j=1}e_i\nabla^A_{e_i}e_j\nabla^A_{e_j}\psi&=&     \sum_{i=j}e_ie_j\nabla^A_{e_i}\nabla^A_{e_j}\psi+\sum_{i,j}e_i     \nabla^A_{\nabla^T_{e_i}e_j}\psi+     \sum_{i<j}e_ie_j\nabla^A_{e_i}\nabla^A_{e_j}\psi+     \sum_{i>j}e_ie_j\nabla^A_{e_i}\nabla^A_{e_j}\psi\\     &=&-\sum_k(\nabla^A_k\nabla^A_k+\nabla^A_{\nabla^T_{e_k}e_k}+\nabla^A_{II(e_j,e_j)}     )\psi+     \sum_{i<j}(e_ie_j)     (\nabla^A_{e_i}\nabla^A_{e_j}     -\nabla^A_{e_i}\nabla^A_{e_i})\psi\\     &=&-\sum_k(\nabla^A_k\nabla^A_k+\nabla^A_{\nabla^T_{e_k}e_k})\psi     +\nabla^A_{\tau_b}\psi+\frac14Scal^T\psi+\frac12F_A\cdot\psi,   \end{eqnarray*}   where $II(e_i,e_i)$ denotes the 2nd fundamental form.   Here to compute the curvature part, we use the classical way, the exactly same way to show the classical Weitzenb\"ock formula. \end{pf}

% \begin{lemma}[Kordyukov, Lejmi and  Weber \cite{KLW}]  Suppose that $(M,F,g)$ satisfies the same hypothesis as the above proposition. Then, the following inequality establishes:   \[\Delta_b(|\psi|^2)\leq2Re((\nabla^A_b)^*\nabla^A_b\psi,\psi),\]   for any pair $(A,\psi)$, where $A$ is a basic \spinc connection and $\psi\in\Gamma_b(S^+)$. \end{lemma}

 For the convenience, we denote
   by $S^T=Scal^T-\delta_b\kappa_b+|\kappa_b|^2$. By the method of \cite[Lemma 16]{KLW}, we have the following lemma.

 \begin{lemma}
   If $(A,\psi)$ is the solution to the basic Seiberg Witten equations \eqref{eqn-SW},  then  we have
   \[|\psi^2|\leq\max_M(-S^T,0).\]
 \end{lemma}

% \begin{pf}   The method is similar to the classical bound formula.   \begin{eqnarray*}     \Delta_b|\psi|^2&\leq&2Re((\nabla^A_b)^*\nabla^A_b\psi,\psi)\\     &=&-\frac12S^T|\psi|^2-Re(F^+_A\psi,\psi)\\     &=&-\frac{S^T}2|\psi|^2-\frac12|\psi|^4.   \end{eqnarray*}   Let $x_0$ be a point where $|\psi|^2$ is maximal. Then, $\tau|\psi|^2=0$, i.e.   $\Delta^{g_Q}|\psi^2|=\Delta_b|\psi|^2$, where   $\Delta^{g_Q}$ is the usual Laplacian associated with the metric   $g_Q$.  Thus, we get the result. \end{pf}

%By using the saluted neighborhood and the partition of unity, we get the basic Sobolev embedding and Sobolev multiplication.\begin{prop}  Setting $0\leq l\leq k$, we have the following statements.  \begin{itemize}    \item[(1)] If $k-\frac{q}{s}>0$ and $k-\frac{q}{s}\geq l-\frac{q}{t}$, then we have        \[L^s_k\times L^t_l\to L^t_l.\]    \item[(2)] If $0\leq l\leq k$ and $l-\frac{q}{t}\leq    k-\frac{q}{s}<l$, then we have    \[L^s_k\hookrightarrow L^t_l.\]    Moreover, for $l<k-\frac{q}{s}$, we have    \[L^s_k\hookrightarrow C^l.\]    \item[(3)] The inclusion    \[L^s_{k+1}\hookrightarrow L^s_k\]    is compact.  \end{itemize}\end{prop}
The arguments of \cite[Lemma 22]{KLW} imply the following proposition.
\begin{prop}
  Let $\mathcal C=L^2_k(\mathcal A_b\times \Gamma_b(S^+))$ for some positive integer $k\geq2$, we have that the quotient  $\mathcal B=\mathcal C/L^2_{k+1}(C^\infty_b(M,S^1))$ is a Hausdorff space.
\end{prop}

 We define the map
 \[\tilde \mu:\mathcal A_b\times(\Gamma_b(S^+)\oplus\Omega^1_b(M)\oplus H^0_b(M) )
 \to \mathcal A_b\times (\Gamma_b(S^-)\oplus \Omega^+_b(M)\oplus\Omega^0_b(M)\oplus \ker(\delta_b+d^+_b)\]
 by
 \[(A,\phi,a,f)\mapsto (A,\Dirac^{A+a}_b\phi,F^+_{A+a}-q(\phi),\delta_ba+f,a_{ha}),\]
 where $\mathcal A_b$ denotes the set of all basic \spinc connections and $a_{ha}$ denotes the image of the $L^2$-projection to the $\ker(d^+_b+\delta_b)$.
 Let $\mathcal G^0_b$ be the fixed-point gauge group.
 Fix a connection $A\in \mathcal A_b$, the subspace $A+\ker(d_b)\subset \mathcal A_b$ is a $\mathcal G^0_b$-invariant subspace.
 We set the following the notations for the later arguments:
 \begin{itemize}
    \item $\tilde{\mathcal A}=(A+\ker(d_b))\times(\Gamma_b(S^+)\oplus\Omega^1_b(M)\oplus H^0_b(M))$,
    \item $\tilde{\mathcal C}=(A+\ker(d_b))\times(\Gamma_b(S^-)\oplus \Omega^+_b(M)\oplus\Omega^0_b(M) \oplus \ker(d^+_b+\delta_b))$
    \item $\mathcal A^b=\tilde{\mathcal A}/\mathcal G^0_b,$
    \item $\mathcal C^b=\tilde{\mathcal C}/\mathcal G^0_b.$
  \end{itemize}

  Since $\tilde \mu$ is $\mathcal G^0_b$--equivariant, we have that  $\mu=\tilde\mu/\mathcal G^0_b:\mathcal A^b\to\mathcal C^b$ is well-defined. The proposition below is  the foliated version of \cite[Proposition 3.1]{BF}

\begin{prop}
  The preimage $\mu^{-1}(B)\subset\mathcal A^b_k$ of a bounded disk bundle $B\subset \mathcal C^b_{k-1}$ is contained in a bounded disk bundle, where $\mathcal A^b_k$ is the $L^2_k$-completion of $\mathcal A^b$ and $\mathcal C^b_{k-1}$ is the $L^2_{k-1}$-completion of $\mathcal C^b$.
\end{prop}

{\bf Remark:}

If $H^1_b(M)\cap H^1(M,\mathbb Z)$ is a lattice of $H^1_b(M)$, we define basic Picard group as $Pic_b(M)=H^1_b(M)/H^1_b(M)\cap H^1(M,\mathbb Z)$.  This  condition is very sensitive to the perturbation of foliation. For example:
   Let $\mathbb Z^2$ be the canonical action on $\mathbb R^n$, for a family of lines with the same slope parameterizing by $\mathbb R$, i.e. $l_\alpha: y=kx+\alpha$ as $\alpha\in\mathbb R$. For a given $k$, one gets a foliation $(\mathbb R^2/\mathbb Z^2,F_k)$, then $H^1_b(M)\cap H^1(M,\mathbb Z)$ is a lattice if and only if $k\in\mathbb Q$.
 Follow the arguments of \cite[Corollary 3.2]{BF}, one establishes the following corollary.

\begin{cor}
Suppose that   $H^1_b(M)\cap H^1(M,\mathbb Z)$ is a lattice of $H^1_b(M)$.
  Then, the monopole map defines an element $[\mu]$ in the stable cohomotopy group,
  \[\pi^0_{S^1,H}(Pic_b(M);\lambda)=\pi^{b^+}(Pic_b(M);Ind(D_b)),\]
  where $H$ is a Sobolev completion of $\Gamma_b(S^-)\oplus\Omega^+_b(X)$, $D_b=\Dirac^A_b\oplus(d^+_b+\delta_b)$ and  $\lambda=Ind(\Dirac)\ominus H^+_b(M)$ is the difference of the complex virtual index bundle of the basic Dirac operator and the trivial bundle $H^+_b(M)$. Here  $S^1$  acts on $Ind(\Dirac)$ as complex multiplication and on $H^+(X)$ trivially.
\end{cor}
\begin{pf}
  Since the operator $\mu=l+c$ satisfies the hypothesis of Theorem \ref{thm-2.6}, it is known that $l+c$ defines a stable cohomotopy Euler class
  $[\mu]\in\pi^0_{S^1,H}(Pic_b(M);\lambda)$.
\end{pf}

Similarly, we have the following foliated version for \cite[Proposition 3.3]{BF}.

\begin{prop}
  Let $(M,F)$ be a closed oriented manifold with codimension $4$  foliation, moreover we assume that $M$ admits a taut bundle-like metric $g$,  $H^1_b(M)\cap H^1(M,\mathbb Z)$ is a lattice of $H^1_b(M)$ and  $b^+_b>b^1_b+1$. By the choice of the homological orientation, i.e. the orientation of $\det(H^1_b(M))\otimes\det(H^+_b(M))$, we have a homomorphism $\mathfrak t:\pi^{b^+}_{S^1,H}(Pic_b(M), Ind(D_b))\to \mathbb Z$, whose value is the Seiberg-Witten invariant.
\end{prop}

At the end of this section, we give a method to construct the examples $(M,F)$ satisfying Assumption \ref{assum-main-1}.
 One way to construct the foliation is by suspension, here we give two   references, see \cite[Chapter 3.8]{Mol} and \cite{Richard}. Let $(Y,g)$ be a closed oriented $4$ Riemannian manifold. Suppose that a compact Lie group $G$ actions on $(Y,g)$ isometrically and preserving the orientation of $Y$, and we have a representation
\[f: \pi_1(X)\to G\]
such that the closure of $Im(f)$ is $G$, where
  $X$ is a closed oriented manifold with fundamental group $\pi_1(X)$.
We set $M=\tilde X\times Y/f$, where $\tilde X$ denotes the universal covering of $X$.
Before preceding, we need the lemma below.

\begin{lemma}[c.f.  {\cite{Richard}}]
  Let $(M,F)$ be defined as above. Then, there is a one-to-one corresponding  $ \Omega^G(Y)\cong \Omega_b(M)$, where $\Omega^G(Y)$ denotes the set of $G$-invariant forms.
\end{lemma}
%\begin{pf}  Setting $\rho:\tilde M=\tilde X\times Y\to M$, we have that $\tilde M$ admits a $\pi_1(X)$ action, defined by $\gamma\cdot(x,y)\mapsto (x\gamma,f(\gamma^{-1})y)$, for any $\gamma\in\pi_1(X)$.   This implies that the following diagram commutes.  \[\xymatrix{\tilde M\ar[r]^{\gamma}\ar[d]&\tilde M\ar[d]\\  M\ar@{=}[r]&M}\]  Since $(M,F)$ lifts a foliation $\tilde F$ of $\tilde M$, we can define the basic forms of $\tilde M$.  Recall that $\Omega_b(M)=\{w\in\Omega^*(M)|~\iota_\xi\omega\equiv0,~L_\xi(\omega)\equiv0,\mbox{ for any }\xi\in\Gamma(F)\}$. We have that the image $\rho^*\Omega_b(M)$ is equal to the $\Omega^{f}_b(\tilde M)$, where  $\Omega^{f}_b(\tilde M)$ denotes the basic forms of $\tilde M$ which are invariant under $f(\pi_1(X))$. Since $\Omega^*_b(\tilde M)=p^*\Omega^*(Y)$, where $p:\tilde M\to Y$, we have that $\rho^*\Omega_b(M)$ is equal to pull-back of $f(\pi_1(X))$-invariant forms of $Y$. By the dense condition, this implies that $\rho^*\Omega_b(M)$ is equal to pull-back of $G$-invariant forms of $Y$.\end{pf}

Since one can find a $G$-invariant volume form over $Y$, lifting back on $M$ we have that $H^3_b(M, F)\neq0$, which implies  that $(M, F)$ admits a taut bundle-like metric.

\begin{lemma}
  Let $(M,F)$ be defined as above. Then,  we have an identification
   \[\pi_0(Map^G(Y,S^1))\cong   H^1(M,\mathbb Z)\cap H^1_b(M).\]
\end{lemma}
\begin{pf}
  It is known that $H^1(M,\mathbb Z)\cong \pi_0(Map(M,S^1))$, which means that  for each element $[w]\in H^1(M,\mathbb Z)$ we have a  representation $u:M\to S^1$ of this homotopy class $[u]_{ht}\in\pi_0(Map(M,S^1))$  satisfying the condition:
   \[ [\frac1{2\pi i}u^{-1}du]=[w].\] Therefore, any element $[w]\in H^1(M,\mathbb Z) \cap H^1_b(M)$ corresponds to a representation $u$ of the homotopy class $[u]_{ht}\in\pi_0( Map(M,S^1))$, such that $[\frac1{2\pi i}u^{-1}du]\in H^1_b(M)$. This implies that there is $f\in i\Omega^0(M)$ such that $u^{-1}du+df\in \Omega^1_b(M)$. Setting $u'=e^{f}u$, we have that $L_\xi u'\equiv0$ for any $\xi\in \Gamma(F)$. This implies that
   \[H^1(M,\mathbb Z)\cap H^1_b(M)\cong\pi_0(Map_b(M,S^1)), \]
   where $Map_b(M,S^1)=\{u|u\in Map(M,S^1),~L_\xi u\equiv0,\mbox{ for any }\xi\in \Gamma(F)\}$. By the above argument, we have that $u$ corresponds to a $G$-invariant function on $Y$, we still use the same notation $u$ to express this $G$-invariant function of $Map(Y,S^1)$. Therefore,
 Thus, a $\mathbb Z$-module subset $\Gamma$ of $H^1_b(M)$ is a lattice of   $H^1_b(M)$ if and only if $rank(\Gamma)=b^G_1(Y)$.
\end{pf}

Suppose there is a $G$-invariant \spinc structure,  Given a $G$-equivariant spinor bundle $S'\to Y,$ we construct a foliated spinor bundle $S=E\times \tilde X/f$.  By \cite{Richard}, it is known that there is  an identification\[\Gamma^G(Y,S')\cong\Gamma_b(M,S).\]
\begin{prop}
   Let  $(M,F)$ be a manifold with foliation constructed  as above and $Y$ admits $G$-equivariant spinor bundle. Suppose  $rank(\pi_0(Map^G(Y,S^1)))=b^G_1(Y)$, where $Map^G(Y,S^1)$ denotes the set of  $G$-invariant $S^1$-valued functions and $b^G_1$ denotes the dimension of the first cohomology for the $G$-invariant deRham complex. Then $(M,F)$ satisfies the Assumption \ref{assum-main-1}.  % We can define the  Suppose that $H^1_b(M)\cap H^1(M,)$
\end{prop}

\begin{pf}
  We need to check that $H^1_b(M)\cap H^1(M,\mathbb Z)$ is a lattice of $H^1_b(M)$. We divide the proof into two parts.
  For $ H^1_b(M)$ part,  we claim that $$H^1_b(M)\cong H^{1,G}_{dR}(Y),$$ where $H^{1,G}_{dR}(Y)$ denotes the first deRham cohomology of $G$-invariant forms. The above lemma  shows that there is a one-to-one corresponding  $ K:\Omega^G(Y)\cong \Omega_b(M)$. By the similar argument, we have that the diagram
  \[\xymatrix{\Omega^{*,G}(Y)\ar[r]^d\ar[d]_K& \Omega^{*+1,G}(Y)\ar[d]^K\\ \Omega^*_b(M)\ar[r]&\Omega^{*+1}_b(M)}\] commutes, which implies that $H^*_b(M)\cong H^{*,G}_{dR}(Y)$. Combining the above lemma, the proposition follows.
\end{pf}

An explicit example is given as follows:
 Let $f:\mathbb Z\to Diff(S^2)$ be a homomorphic map, generated by the
rotating $S^2$ around the $z$-axis at $\alpha$ degree, where $\alpha$ is not a
rational number.
  We set $$M_0=\mathbb R\times S^2/\sim_f,$$ where the equivalent relation $\sim _f$ is defined  as follows: we say two points $(x_1,y_1)$ and $(x_2,y_2)$ are equivalent
   if and only if there exists $n\in\mathbb Z$ such that $x_1+n=x_2$ and $f(-n)y_1=y_2$, we   denote this equivalent class by $[x,y]$.
  Then the projection to the first component becomes a foliation, whose leaf
at $[x,y_0]$ consists of all such points $(t,y)$, where $t\in[0,1]$ and $y=f(n)y_0$ for all $n\in\mathbb Z$.
  Equip $S^2$ with the canonical metric, we have that $(M_0,\mathcal F)$ is a Riemannian foliation, furthermore the leaves over generic point are non-compact. Let $\Sigma_g$ be the closed Riemannian surface with genus $g$ and $M=M\times \Sigma_g$.  One has that $(M,\mathcal F)$ is a codimension-$4$ Riemannian foliation, $H^1_b(M)\cap H^1(M,\mathbb Z)$ is
a lattice in $H^1_b(M)$ and $TM/T\mathcal F$ admits a transversal spin structure.

\section{Proof of Theorem \ref{thm-main-estimate}}

In this section, we assume that a codimension-$4$ foliation of $(M,F)$ admits a foliated spin structure. Fixing a bundle like metric $g_Q$, one has the associated transversal Dirac operator, $\Dirac_b$. We define  the operators $l,~c:i\Omega^1_b\times\Gamma_b(S^+)\to i\Omega^+_b\times\Gamma_b(S^-)\times i\Omega^0$ by
\[\begin{array}{l}
  l(a,\Phi)=(d^+_ba,\Dirac_b\Phi,\delta_ba),\\
  c(a,\Phi)=(-\sigma(\Phi),\frac a2\cdot\Phi,0).\\
\end{array}\]
 Let $\mathbb H$ be the quaternion numbers and $Sp(1)$ be the group of the quaternion numbers with norm $1$, let $Pin(2)$ be the normalizer of $S^1$ in $Sp(1)$, i.e. $Pin(2)=S^1\cup jS^1\subset Sp(1)\subset\mathbb H$, it acts on $\Gamma_b(S^\pm)$ by the right multiplication
\begin{equation}
  \begin{array}{rcl}
  Pin(2)\times\Gamma_b(S^\pm)&\longrightarrow&\Gamma_b(S^\pm),\\\
  (g,\Phi)&\longmapsto&\Phi\cdot g^{-1},\\
\end{array}\label{formula-Pin2-action}
\end{equation}
where $g$ is any element of $Pin(2)$ and $\Phi$ is any element of $\Gamma_b(S^\pm)$; it acts on the spaces $\Omega^1_b$ and $\Omega^+_b$ as follows:
 the $S^1$ component acts trivially and $j$ acts as $-1$-multiplication, i.e. $Pin(2)$ acts via the canonical projection $Pin(2)/S^1\cong \mathbb Z_2$ as $\pm1$ multiplication.
Here $l$ is a variant of the operator $D$ and $c$ is a variant of  the operator $Q$ in Furuta's paper \cite{Fu}. In particular, when the foliation is trivial, i.e. $0$-dimension, we have that  $l=D$ and $c=Q$.
%\begin{lemma}[Furuta \cite{Fu}]  $l,~c$ are $Pin(2)$-equivariant.\end{lemma}
By the straightforward calculation(see  \cite{Fu}), we have that $l,~c$ are $Pin(2)$-equivariant.
Let $V_b=L^2_4(\Omega^1_b\oplus \Gamma_b(S^+)),~W_b=L^2_3(i\Omega^+_b\oplus \Gamma_b(S^-)\oplus i\Omega^0_b)$. By Sobolev multiplication theorem, we have $L^2_3\otimes L^2_3\to L^2_3$. %Let $D^*_b$ be adjoint of $D_b$ under such a metric, for the operator $D$ defined in the previous subsection.
We define the norms of $V_b$ and $W_b$ respectively,
\[\|v\|^2_{V_b}=\int_X\left(|(l^*l)^2v|^2+|v|^2\right),~~
\|w\|^2_{W_b}=\int_X\left(|(l^*l)^{\frac32}w|^2+|w|^2\right).\] We define the $Pin(2)$-action on $V_b$ and $W_b$ as follows: on the spinor section component, $Pin(2)$ acts as \eqref{formula-Pin2-action}, and  on the forms component, it acts as $\pm1$ multiplication.
We have that
\[\|s\|_{L^2_3}\mbox{ is equivalent to }\left(\|s\|_{L^2}+\|(l^*l)^{\frac32}s\|_{L^2},\right)\]
for any $s\in L^2_3.$

%\begin{lemma}[Furuta \cite{Fu}]  $\|\cdot\|_V,~\|\cdot\|_W$ are $Pin(2)$ invariant.\end{lemma}

%Since $(l+c)^{-1}(0)/S^1$ is the normal moduli space of Seiberg-Witten equation, we have the following lemmas.

%\begin{lemma}[Furuta \cite{Fu}]\label{lemma-compact-zero-locus}  Zero locus of  map $l+c:V_b\to W_b$ is compact.\end{lemma}

\noindent Since we have that  $l^*l,~ll^*\geq0$, for each nonnegative real number $\lambda$, we denote $V_{b,\lambda}(W_{b,\lambda})$ the eigen-space consisting of eigenvector with eigenvalue  $\leq\lambda$ of $ll^*$($l^*l$). Let $p_\lambda:W_b\to W_{b,\lambda}$ be  $L^2$-projection.

\noindent
We set the operator
\[l_\lambda+c_\lambda:V_\lambda\to W_\lambda,\]
where $l_\lambda$ denotes the restriction of $l$ and $c_\lambda=pr_\lambda\comp c$.

\noindent
We define operators,
  \[D_1=\Dirac_b:L^2_4(S^+)\to L^2_3(S^-),~
  D_2=d^+_b+\delta_b:L^2_4(i\Omega^1)\to L^2_3(i\Omega^+\oplus i\Omega^0).\]
  Let $V^1_\lambda$ be the subspace consisting of eigenvectors  with corresponding eigenvalues $\leq\lambda$ of $D^*_1D_1$ and $W^1_\lambda$  be the subspace consisting of eigenvectors  with corresponding eigenvalues $\leq\lambda$ of $D_1 D^*_1$. Similarly, we set the spaces $ V^2_\lambda$ and $W^2_\lambda$ as the subspaces consisting of eigenvectors  with corresponding eigenvalues $\leq\lambda$ of $ D^*_2 D_2$  and $D_2 D^*_2$ respectively.
\noindent
  We regard $V_\lambda=V^1_\lambda\oplus V^2_\lambda,~W_\lambda=W^1_\lambda\oplus W^2_\lambda$ as $Pin(2)$ action decomposition.  Since $D^*_1D_1=\Dirac^*_b\Dirac_b$ commutes with $Pin(2)$, we have that $V^1_\lambda$ is a $\mathbb H$--linear space, i.e. $V^1_\lambda=\mathbb H^{m'}$ for some $m'\in\mathbb N$. Similarly, $W^1_\lambda=\mathbb H^m$ for some $m\in\mathbb N$.
  By straightforward calculation,  we have
  \[\dim_{\mathbb R}V^1_\lambda-\dim_{\mathbb R}W^1_\lambda=
  \dim_{\mathbb R}V^1_0-\dim_{\mathbb R}W^1_0=\dim_{\mathbb R}\ker(D^*_1D_1)
  -\dim_{\mathbb R}\ker(D_1D^*_1).\]
  By the relations $\ker(D^*_1D_1)=\ker(D_1)$ and $\ker(D_1D^*_1)=\ker(D^*_1)$, we have
  \[\dim_{\mathbb R}\ker(D^*_1D_1)-\dim_{\mathbb R}\ker(D_1D^*_1)=Ind_{\mathbb R}(D_1).\]
   We get that $Ind_{\mathbb R}(D_1)=\dim_{\mathbb R}V^1_\lambda-\dim_{\mathbb R}W^1_\lambda=4m'-4m$, we  denote by $4k=Ind_{\mathbb R}(D_1)$.
  Now, we consider $V^2_\lambda,~W^2_\lambda$. Since $Pin(2)$ acts on $V^2_\lambda,~W^2_\lambda$ as $\mathbb Z_2$,  we   write $V^2_\lambda=\mathbb R^n$ and $W^2_\lambda=\mathbb R^{n'}$. By the definition $D_2=d^+_b+\delta_b:i\Omega^1_b\to i\Omega^+_b\oplus i\Omega^0_b$, we have
  \[Ind_{\mathbb R}(D_2)=n-n',\]
  and by \cite[Proposition 3]{KLW}, we have that  $Ind_{\mathbb R}(D_2)=-(b^0_b-b^1_b+b^+_b)=-(b^+_b+1)$,
   $n'=b^+_b+n+1$, i.e. $W^2_\lambda=\mathbb R^{b^+_b+n+1}$.
  Therefore by the above lemma, we have  a $Pin(2)$-equivariant map
  $l_\lambda+c_\lambda:\mathbb H^{r+k}\oplus\mathbb R^n\to  \mathbb H^{r}\oplus\mathbb R^{n+b^+_b+1}$.

On has the following lemma.

\begin{lemma}[Furuta {\cite[Lemma 3.4]{Fu}}]\label{lemma-3.4}
 For  large enough $\lambda$, the operator $l_\lambda+c_\lambda$ has no zeros on the sphere with radius $R$
centered in $0$.
\end{lemma}
The above arguments also implies  the following lemma.

\begin{lemma} Let $l_\lambda+c_\lambda$ be defined as above. Then,
for  large enough $\lambda$, the operator
  $l_\lambda+c_\lambda:V_\lambda\to W_\lambda$ is a $Pin(2)$-equivariant map, and it defines a smooth map $(BV_\lambda,SV_\lambda)\to (B \bar W_\lambda,S\bar W_\lambda)$, where $BV_\lambda$($BW_\lambda$) denotes  the unit ball of $V_\lambda$($\bar W_\lambda$) and $SV_\lambda$($S\bar W_\lambda$) denotes the unit sphere of $V_\lambda$($\bar W_\lambda$)
\end{lemma}
\begin{pf}
  By Lemma \ref{lemma-3.4}, we may assume that $l_\lambda+c_\lambda$ has no zero point on $B_R(V_\lambda)\setminus B_\delta(V_\lambda)$ for some $\delta<R$, where $B_R(V_\lambda)$ denotes the ball with radius $R$ centered at the origin. Since $l_\lambda+c_\lambda$ is smooth and $BV_\lambda\cup SV_\lambda$ is a compact subset of $V_\lambda$, there is a positive number $C>0$, such that $Im((l_\lambda+c_\lambda)|_{BV_\lambda})\subset B_C(\bar W_\lambda)$. We consider the map
  \[F:(BV_\lambda,SV_\lambda)\to (B\bar W_\lambda,S\bar W_\lambda)\]
  defined by
  $F(\gamma)=\rho(\|\gamma\|_W) \frac{(l_\lambda+c_\lambda)(\gamma)}C+(1-\rho(\|\gamma\|))
  \frac{(l_\lambda+c_\lambda)(R\cdot\gamma)}{\|(l_\lambda+c_\lambda)(R\cdot \gamma)\|_W}$, where $\rho:\mathbb R\to\mathbb R$ is a cutoff function, such that $\rho(t)=1$ for $|t|\leq \delta/2$ and $\rho(t)=0$ for $|t|\geq \delta'$ for some $\delta'\in (\delta, R)$. Clearly, $F$ is $Pin(2)$-equivariant.
\end{pf}

We show that the image of the operator $l_\lambda+c_\lambda$ is contained in a codimension $1$ subspace $\bar W_\lambda$ of $W_\lambda$. The reason is not different   to \cite[Remark before Theorem 4.2]{Fu}. Here we give a brief explanation. It is known that $W_b$ contains a parallel $0$-form $s_0$, i.e. the constant function on $M$. Since $\int_M\delta_ba=0$ for any $a\in i\Omega^1_b$, it is clear that the this section is contained in the kernel of the operator $D^*_2$(the formal adjoint operator to $D_2$). The image of $D_2$ is contained in the  $L^2$-orthogonal complement to $s_0$ in $W_b$. The image of $c_\lambda$ is also in the  $L^2$-orthogonal complement to $s_0$ in $W_b$.  From the construction of the finite dimensional approximation, it is known that the image of $l_\lambda+c_\lambda$ is still contained in the subspace $\bar W_\lambda=W_\lambda\cap s^\perp_0\cong \mathbb H^{r}\oplus\mathbb R^{n+b^+_b}$ , where $s^\perp_0$ denotes the $L^2$-orthogonal complement to $s_0$ in $W_b$.

%\begin{pf}  As  $\|v\|_V=R$, by Lemma \ref{lemma-finite-approximation-2},  \[\epsilon\leq\|(l+c)v\|_W\leq\|(l+p_\lambda c)v\|_W+\|p^\lambda cv\|_W,\]  combining with Lemma \ref{lemma-finite-approximation-3}, we have  \[0<\epsilon-\|p^\lambda c(v)\|_W\leq\|(l+p_\lambda c)v\|_W.\]  Thus, $l+p_\lambda c:V\to W$ has no zero point on the ball with radius $R$, after taking the restriction, the lemma is proved.\end{pf}

%\begin{lemma}\label{lemma-approximation-decomposition}

%\end{lemma}

Combining the above lemma and proposition, we give a proof of Theorem \ref{thm-main-estimate}.

\begin{pf}( of Theorem \ref{thm-main-estimate})
We may assume that $Ind(\Dirac_b)\geq0$, otherwise the theorem is trivial.
  By the previous arguments, we have  a $Pin(2)$-equivariant map
  $l_\lambda+c_\lambda:\mathbb H^{r+k}\oplus\mathbb R^n\to  \mathbb H^{r}\oplus\mathbb R^{n+b^+_b}$.% is, and it defines a smooth map $(B(\mathbb H^{r+k}\oplus\mathbb R^n),S(\mathbb H^{r+k}\oplus\mathbb R^n))\to (B (\mathbb H^{r }\oplus\mathbb R^{n+b^+_b}),S(\mathbb H^{r }\oplus\mathbb R^{n+b^+_b}))$.
   There is an induced map between
the balls in the corresponding complexified representations and the
proposition below shows that such a map can only exist if either $k = 0$ or $b^+_b\geq 2k+1$.
\end{pf}

\begin{prop}[Furuta {\cite[Proposition 5.1]{Fu}}]
 Let   \[V =\mathbb H^{k+m}\oplus\mathbb C^n,~W_\lambda=\mathbb H^m\oplus\mathbb C^{b +n},\]
  and suppose that there is a $Pin(2)$-equivariant smooth map $f:(BV,SV)\to (BW,SW)$, where $BV$($BW$) denotes  the unit ball of $V$($W$) and $SV$($SW$) denotes the unit sphere of $V$($W$). Then, either $k=0$ or $b\geq 2k+1$.
\end{prop}

 %%%%%%%%%%%%%

Graduate School of Mathematical Sciences, The University of Tokyo, 3-8-1 Komaba, Meguro-ku, Tokyo 153-8914, Japan. E-mail: dexielin@ms.u-tokyo.ac.jp

\end{document}